\documentclass[a4paper,12pt,leqno]{article}
\usepackage{amsmath,amsfonts,amsthm, amssymb}
\usepackage[OT4]{fontenc}
\usepackage{graphicx}
\usepackage{longtable}
\usepackage{setspace}
\usepackage[left=3cm,top=2.5cm,right=3cm]{geometry}

\long\def\symbolfootnote[#1]#2{\begingroup%
\def\thefootnote{\fnsymbol{footnote}}\footnote[#1]{#2}\endgroup}

\newtheorem{theorem}{Theorem}[section]
\newtheorem{lemma}[theorem]{Lemma}
\newtheorem{proposition}[theorem]{Proposition}
\newtheorem{cor}[theorem]{Corollary}
\newtheorem{wtw}{Theorem}

\theoremstyle{definition}
\newtheorem{rem}[theorem]{Remark}
\newtheorem{obs}[theorem]{Observation}
\newtheorem{defin}[theorem]{Definition}

\renewcommand{\proof}{\medskip\par\noindent\textbf{Proof.} \ignorespaces}

\renewcommand{\qed}{\quad\hskip0pt\null\hfill$\square$\par}

\newcommand{\R}{\mathbb{R}}

\newcommand{\p}{\partial}

\def\.{\hskip.06cm}
\def\ts{\hskip.03cm}

\begin{document}

\begin{center}
\large\bfseries
Acute triangulations of polyhedra and $\R^n$
\end{center}

\begin{center}\bf
Eryk Kopczy\'nski,
Igor Pak\symbolfootnote[1]{Partially supported by the
University of Minnesota, Universit\'e Paul Sabatier and the NSA.} \&
Piotr Przytycki\symbolfootnote[2]{Partially supported by MNiSW grant
N201 012 32/0718, the Foundation for Polish Science, and ANR grant ZR58.}
\end{center}

\begin{abstract}
\noindent We study the problem of \emph{acute triangulations} of
convex polyhedra and the space~$\R^n$. Here an acute triangulation
is a triangulation into simplices whose dihedral angles are acute.
We prove that acute triangulations of the $n$--cube do not exist for
$n\ge 4$. Further, we prove that acute triangulations of the
space~$\R^n$ do not exist for $n\ge 5$. In the opposite direction,
in~$\R^3$, we present a construction of an acute triangulation of the cube,
the regular octahedron and a non-trivial acute triangulation of
the regular tetrahedron. We also prove nonexistence of an
acute triangulation of~$\R^4$ if all dihedral angles are bounded
away from~$\pi/2$.
\end{abstract}



\section{Introduction}

The subject of acute triangulations is an important area of Discrete
and Computational Geometry, with a number of connections to other areas
and some real world applications. Until recently, most results dealt with
the $2$--dimensional case, where the problem has been largely resolved.
In the last few years, several papers~\cite{ESU,Kri,VH} broke the dimension
barrier in both positive and negative direction (see below). In this
paper we continue this exploration, nearly completely (negatively)
resolving the problem in dimension~4 and higher, and making further
advancement in dimension~3.

The problem of finding acute triangulations has a long history in
classical geometry, and is elegantly surveyed in~\cite{BS}, which
argues that it goes back to Aristotle. In recent decades, it was
further motivated by the \emph{finite element method} which requires
``good'' meshes (triangulations of surfaces) for the numerical
algorithms to run. Although the requirements for meshes largely depend
on the algorithm, the sharp angle conditions seem to be a common feature,
and especially important in this context. We refer to~\cite{Str} for
the introduction to the subject, and to~\cite{Sw} for the state of art.

Another motivation comes from the recreational literature,
where the subject of dissections has been popular in general
(see~\cite{Lin}), and of acute triangulations in particular~\cite{CL,Man}.
In this context, the problem of acute triangulations of a square, cube,
and hypercubes seem to be of special interest~\cite{Epp}.

\medskip

An \emph{acute triangulation} is a dissection into \emph{acute}
simplices (i.e.\ with acute dihedral angles) which form a simplicial
complex, so e.g. in the plane, a vertex of one simplex cannot lie in
the interior of an edge of another. (See Figure~1, where on the left
we have a dissection of the square which forms a simplicial complex,
and on the right we have a dissection which does not.)

\begin{center}
\includegraphics{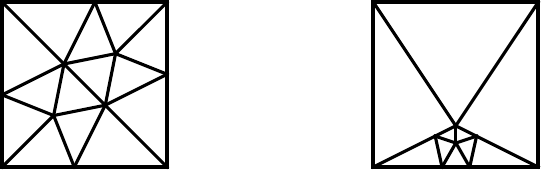}
\centerline{Fig. 1}
\end{center}

In one of the first papers on the subject, Burago and Zalgaller
proved in~\cite{BZ} that any non-convex polygon (possibly, with
holes) has an acute triangulation. Unfortunately, their argument was
largely forgotten as it was inexplicit and did not give a bound on
the number of triangles required. In a long series of
papers~\cite{BE,BGR,BMR,M,Sar,Y} first polynomial, and then linear
bounds were obtained for non-obtuse, and, eventually, for acute
triangulations. We refer to~\cite{Zam} for the historical outline, a
short survey, and further references.

In higher dimensions, several results have been recently obtained.
First, Eppstein, Sullivan and \"{U}ng\"{o}r \cite{ESU} showed that
the space~$\R^3$ can be triangulated into acute tetrahedra, by
adopting a classical tiling construction due to Sommerville. Then,
K{\v{r}}{\'{\i}}{\v{z}}ek~\cite{Kri} showed that no vertex in~$\R^n$
for~$n\ge 5$ can be surrounded by a finite number of acute
simplices\footnote{There is a crucial error in this proof. We
refer to Subsection~\ref{Krizek} for the details.} and conjectured
that the space~$\R^4$ also cannot be triangulated into acute
tetrahedra. Finally, and most recently, VanderZee, Hirani,
Zharnitsky and Guoy~\cite{VH} used an advanced numerical simulation
technique to find an acute triangulation for the (usual) cube
in~$\R^3$. Their construction is independent of ours and uses fewer
tetrahedra.

\medskip

In this paper we prove several results in higher dimensions.

\begin{wtw}[Theorems \ref{dim3} and \ref{dim3tetrahedra}(i)]
\label{1}
There exists an acute triangulation of the cube,
the regular octahedron, and a non-trivial
acute triangulation of the regular tetrahedron.
\end{wtw}

Roughly, we first triangulate the cube into a regular $3$--simplex
and four standard $3$--simplices (and the octahedron into eight
standard ones). We then subdivide each of these $3$--simplices into
$543$ pieces, to obtain combinatorially what we call the
\emph{special subdivision} (based on the $600$--cell, see
Section~\ref{construction}). This approach was used previously by
Przytycki and \'Swi\k{a}tkowski in~\cite{PS} to construct the so
called \emph{flag-no-square} subdivisions in dimension~$3$ (see
Definition~\ref{fns}). Let us mention here that this ``curvature''
condition was surveyed in the appendix of~\cite{PS}, and that it was
used originally to construct Gromov hyperbolic groups with
prescribed boundaries. Let us repeat that the case of the cube
in~$\R^3$ was independently resolved in~\cite{VH}.

\medskip

In the opposite direction, we prove the following result:

\begin{wtw}[Corollary \ref{4-cube}]
\label{2}
There is no periodic acute
triangulation of the space~$\R^4$.
In particular, there is no acute triangulation of the $4$--cube.
\end{wtw}

The first assertion of Theorem~\ref{2} implies
the second one, as $4$--cubes tile the space
(see Section~\ref{Rich triangulations of 4--manifolds}).
A short combinatorial proof of Theorem \ref{2} is based
on the generalized Dehn--Sommerville equations. This method
also gives new results on flag-no-square triangulations
(see Section~\ref{Rich triangulations of 4--manifolds}).
Moreover, it allows to complete the acute triangulations
picture with the following.

\begin{wtw}[Corollary \ref{higher dim}, {\cite[Theorem 6.2]{Kri}}]
\label{higher dim introduction}
There is no triangulation of a polyhedron in $\R^n$, for $n\geq 5$, which
contains an interior vertex such that all dihedral angles adjacent to it
are acute.
\end{wtw}

In particular, there is no acute triangulation of $\R^n$ and the
$n$--cube for $n\geq 5$.

\medskip

Finally, we prove the following most general result:

\begin{wtw}[Corollary \ref{away from 90}]
\label{3}
For every $\varepsilon > 0$, there is no triangulation of the
space~$\R^4$ into simplices with dihedral angles
less than $\frac{\pi}{2} - \varepsilon$.
\end{wtw}

The proof of Theorem~\ref{3} relies on the generalized Dehn--Sommerville equations
and on the relations between isoperimetric inequalities and parabolicity of
infinite graphs.

\medskip

The paper is structured as follows. In Section~\ref{construction} we
study acute triangulations in~$\R^3$ and prove Theorem~\ref{1}. In a
short Section~\ref{section D-S} we recall the generalized
Dehn--Sommerville equations. Then, in Section~\ref{Rich
triangulations of 4--manifolds}, we study their consequences for
\emph{rich} triangulations (combinatorial consequence of both acute
and flag-no-square, see Definition \ref{rich}), and prove
Theorem~\ref{2} and Theorem \ref{higher dim introduction}. We then
switch our attention to Theorem~\ref{3} in Section~\ref{Acute-angled
triangulations of R4}.

\medskip \noindent \textbf{Convention.} In the entire article we
adopt a convention that simplicial complexes and triangulations of
(homology) manifolds are not allowed to have edges connecting a
vertex to itself, and they are also not allowed to have multiple
simplices spanned on the same set of vertices.

\bigskip
\noindent
\textbf{Acknowledgements.} This work was initiated at
Universit\'e Paul Sabatier in Toulouse,
where the second and the third author were visiting.
We thank the university and Jean-Marc Schlenker for their hospitality.

We thank Marc Bourdon, who suggested the final argument of
Section~\ref{Acute-angled triangulations of R4}. We are grateful to
Anil Hirani for telling us about~\cite{VH} and explaining the ideas
behind this work, and to Evan VanderZee for giving us his numerical
estimates. We are thankful to Michal K{\v{r}}{\'{\i}}{\v{z}}ek for
confirming the crucial error in his paper~\cite{Kri} and telling us
about a forthcoming correction. We are grateful to Jon McCammond for
pointing out an error in our earlier version of the proof of
Proposition~\ref{isoperimetry}. 
We also thank Itai Benjamini, Isabella Novik, Vic Reiner, Egon
Schulte and Alex Vladimirsky for help with the references.

\section{Acute triangulations of the $3$--cube and the octahedron}
\label{construction}
In this section we describe acute triangulations of the $3$--cube and
the octahedron. The starting point is the following observation:

\begin{obs}
\label{five}
The link of an interior edge of an acute triangulation of a~polyhedron in $\R^3$
is a~simplicial loop of length at least $5$.
\end{obs}

In view of this observation, let us make the following definition.

\begin{defin}
\label{rich} A triangulation of an $n$--dimensional homology 
manifold is \emph{rich} if the links of all interior
$(n-2)$--simplices are loops of length at least $5$.
\end{defin}

Note that being rich is a~purely combinatorial (i.e.\ non-metric)
property. Observation~\ref{five} states that an acute triangulation
of the~$3$--cube (or a regular octahedron) must be rich. We compare
this definition with the following notion:

\begin{defin}
\label{fns} A simplicial complex (or a~triangulation) is called
\emph{flag-no-square}, if it is \emph{flag} (i.e.\ each set of
vertices pairwise connected by edges spans a~simplex) and each
simplicial loop of length four \emph{has a~diagonal} (i.e.\ a~pair
of opposite vertices of the loop spans an edge).
\end{defin}

\begin{rem}
\label{fns->rich}
Every flag-no-square triangulation of a homology manifold is rich. 
\end{rem}

\medskip
Przytycki--\'Swi\k{a}tkowski \cite[Corollary 2.14]{PS} proved that every
$3$--dimensional polyhedral complex admits
a~flag-no-square subdivision. We recall this construction, since we also
use it to subdivide the $3$--cube and the octahedron.

\medskip
Consider the $600$--cell, the convex regular $4$--polytope with Schl\"afli
symbol $\{3; 3; 5\}$ (see, e.g., \cite{C}).
Denote by $X_{600}$ the boundary of the $600$--cell, a $3$--dimensional simplicial
polyhedron homeomorphic to the $3$--dimensional sphere. It consists of~$600$
$3$--simplices\footnote{To streamline and simplify the presentation,
we refer to triangles as 2--simplices, to tetrahedra as 3--simplices, etc.}
and has $120$ vertices. Its vertex links are icosahedra and its edge links are
pentagons. We first focus on the combinatorial simplicial structure of~$X_{600}$.
Denote by $X_{543}$ the subcomplex of $X_{600}$ which we obtain by removing from $X_{600}$
the interiors of all simplices intersecting a~fixed $3$--simplex. (The number $543$ in the
subscript refers to the number of $3$--simplices in $X_{543}$.)

\begin{lemma}[{\cite[Lemmas 2.5 and 2.7]{PS}}]
\label{X543good}
\item[(1)] $X_{543}$ is topologically a $3$--ball. It is flag-no-square.
\item[(2)] Its boundary is a~$2$--sphere simplicially isomorphic to the simplicial complex which we obtain from the
boundary of a~$3$--simplex by subdividing each face as in Figure 2.
\end{lemma}

\begin{center}
\includegraphics{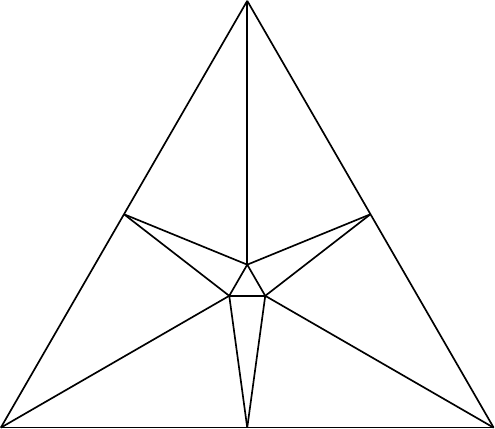}
\centerline{Fig. 2}
\end{center}

We recall the following definition:
\begin{defin}[{\cite[Definitions 2.2 and 2.8]{PS}}]
Given a simplicial complex of dimension at most $3$, its \emph{special subdivision}
is the simplicial complex obtained by:
\item[(i)] subdividing each edge into two (by adding an extra vertex
in the interior of the edge),
\item[(ii)] subdividing each $2$--simplex as in Figure $2$,
\item[(iii)] subdividing each $3$--simplex so that it becomes
isomorphic to $X_{543}$.
\end{defin}

\begin{center}
\includegraphics{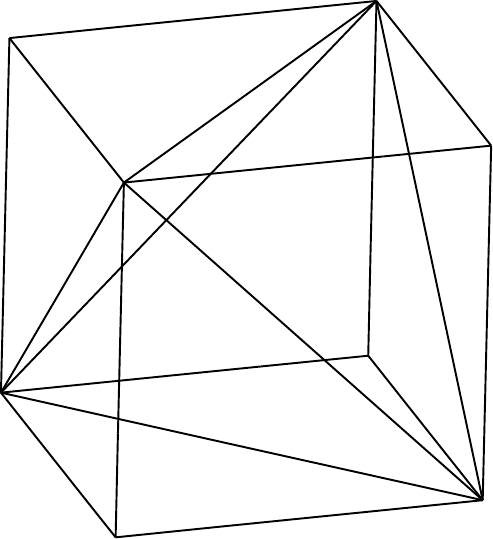}
\centerline{Fig. 3}
\end{center}

We are ready to describe the combinatorial structure of our
triangulation of the $3$--cube and the octahedron. Assume that the
cube lies in $\R^3$ with the vertices at points $(\pm 1, \pm 1, \pm
1)$. Consider the triangulation~$W$ of the cube into five
$3$--simplices so that one of them (denote it by $T_0$) has vertices
$(1,1,1), (-1,-1,1),(-1,1,-1)$ and $(1,-1,-1)$, while the remaining
four $3$--simplices (denote them $T_1,\ldots, T_4$) are the
components of the complement to~$T_0$ in the cube (see Figure $3$).
Note that $T_1,\ldots,T_4$ are congruent (equal up to a rigid
motion); we call such $3$--simplices \emph{standard}
(see e.g.~\cite{Pak}).\footnote{This
tetrahedron is also called the \emph{cube-corner}.} Let $W^*$ be
the special subdivision of~$W$ defined as above.

Similarly, let~$Y$ be the triangulation of the octahedron into eight
standard $3$--simplices obtained as cones from the center over the faces.
Let $Y^*$ be the special subdivision of~$Y$.

By~\cite[Proposition 2.13]{PS}, subdivisions $W^*$ and $Y^*$ are
both flag-no-square. Thus, they are rich and have a potential of
giving an acute realization. This is true indeed, and the main
result of this section is the following theorem:

\begin{theorem}[part of Theorem \ref{1}]
\label{dim3}
\item[(1)] There is an acute triangulation of the $3$--cube, which is combinatorially equivalent to $W^*$.
\item[(2)] There is an acute triangulation of the octahedron, which is combinatorially equivalent to $Y^*$.
\end{theorem}

In fact, we provide acute triangulations of all $3$--simplices of
$W$ and $Y$, combinatorially equivalent to $X_{543}$, and matching
on common part of the boundary. In other words, we prove the
following intermediate result:

\begin{theorem}[part of Theorem \ref{1}]
\label{dim3tetrahedra}
There is a~(non-trivial) acute triangulation, combinatorially equivalent to $X_{543}$, of
\item{(i)} the~regular $3$--simplex,
\item{(ii)} the standard $3$--simplex.
\end{theorem}

Below we describe the construction for the $3$--cube. At some points we use a computer program. We provide the
exact position of all vertices of both triangulations from Theorem \ref{dim3tetrahedra} in the appendix.
There are three steps of the construction. First, we construct an acute triangulation of $T_0$. Then
we ``flatten" it to obtain an acute triangulation of $T_1$. Then we construct another acute
triangulation of $T_0$ so that it matches the one of $T_1$ on the common part of the boundary.

\medskip\par\noindent\textbf{Step 1.}\ignorespaces
\ Note that the vertices of the $600$--cell, whose boundary we
called $X_{600}$, lie on a~sphere in $\R^4$. Moreover all
$3$--simplices in this realization of $X_{600}$ are regular, hence
acute. Let now $\widetilde{X}_{543}$ be the realization of $X_{543}$
in $\R^3$, whose vertices are obtained by stereographic projection
of the $\R^4$ realization. We choose the center of the projection to
be the~center of the (spherical)~$3$--simplex in $X_{600}$ disjoint
from $X_{543}$. It turns out that this mapping does not disturb the
angles significantly.

We move the vertices of $\partial \widetilde{X}_{543}$ radially so
that they arrange on the boundary of a regular $3$--simplex which we
identify with $T_0$. If we scale the size of $\partial T_0$
correctly, this triangulation of $T_0$ is already acute, i.e.\ it
satisfies Theorem \ref{dim3tetrahedra}(i). However, it is not the
one listed in the appendix, we will modify it later in Step~3 (see
also Figure~4).

\begin{center}
\includegraphics{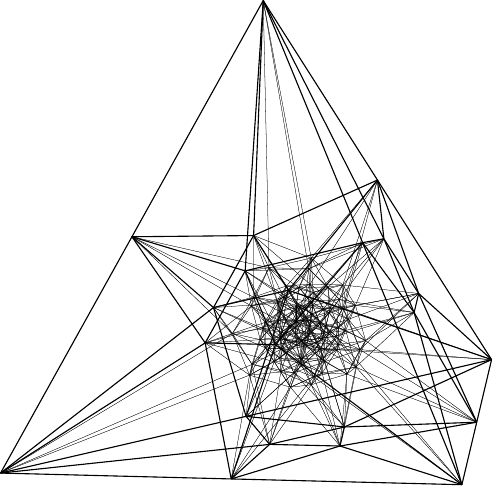}
\centerline{Fig. 4}
\end{center}

\medskip\par\noindent\textbf{Step 2.}\ignorespaces
\ The subdivision of the standard $3$--simplex, say $T_1$, is more difficult.
Our computer program uses the following algorithm to find the position of the vertices.
We ``flatten'' the acute triangulation of $T_0$ obtained in Step 1 in order to obtain
an acute triangulation of $T_1$. We gradually move one of the boundary vertices
(marked~$A$ on Figure~5)
towards the center, keeping the three vertices marked~$C$ in the points
where the circles inscribed into triangles $AB_iB_j$ meet the edges $AB_i$.
We also keep the nine points marked~$D$ on their faces, and scale and
translate together all the interior vertices.

\begin{center}
\includegraphics{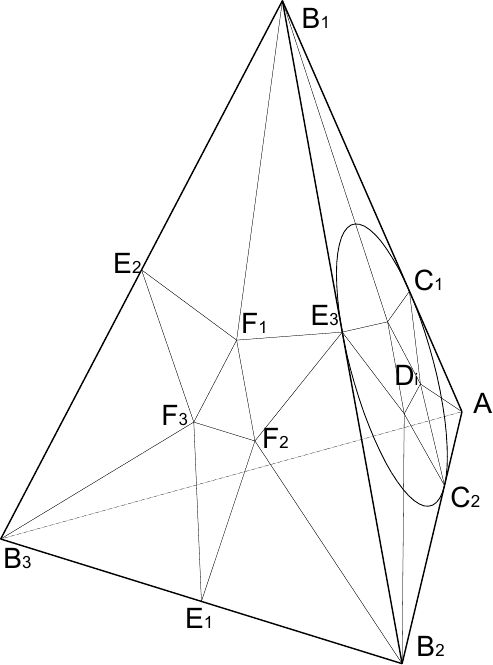}
\centerline{Fig. 5}
\end{center}

Whenever some angle stops being acute during this operation, we
suspend the flattening process to correct the angles. This is done
by slightly moving the responsible vertices so that the angle
becomes smaller. Vertices are moved only in a way that does not
disturb the combinatorial structure, i.e.\ points $A, B_i, C_i, E_i$
are not moved at all; movement of $D_i$ and $F_i$ is restricted to
their faces; all the interior vertices except the two outermost
layers (of $12$ and~$16$ vertices, respectively) are moved together
so that the structure is not disrupted.

When all the angles are corrected, we resume
the flattening, until we obtain the standard $3$--simplex $T_1$.
This completes the description of the triangulation
in Theorem~\ref{dim3tetrahedra}, part~(ii) (see also Figure~6 and the appendix).

\begin{center}
\includegraphics{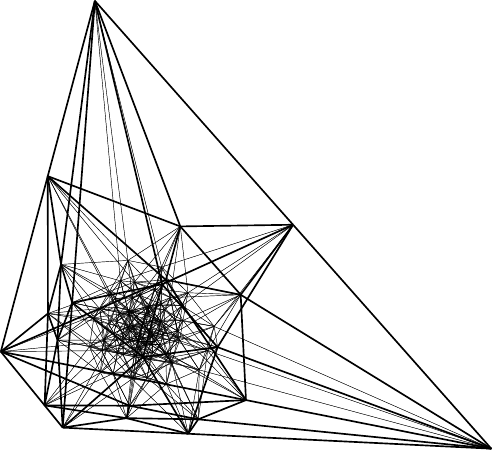}
\centerline{Fig. 6}
\end{center}

\medskip\par\noindent\textbf{Step 3.}\ignorespaces
\ The position of the vertices $F_i$ on the equilateral face of $T_1$ is now different from
their position on the face of $T_0$, because we had to move $F_i$ during the correcting process in Step 2.
So in the triangulation of $T_0$ constructed in Step 1 we move all $12$ vertices corresponding to $F_i$ to the position
matching with the standard $T_i$. It turns out that it is then enough to scale the interior structure to obtain an acute
triangulation (see the appendix).

Now we attach acute triangulations of all $T_i$, constructed in
Steps~2 and~3, to obtain a~triangulation satisfying conditions of
Theorem~\ref{dim3}, part~(1) (see Figure~$7$ and Remark~\ref{Eryk}).

\begin{center}
\includegraphics{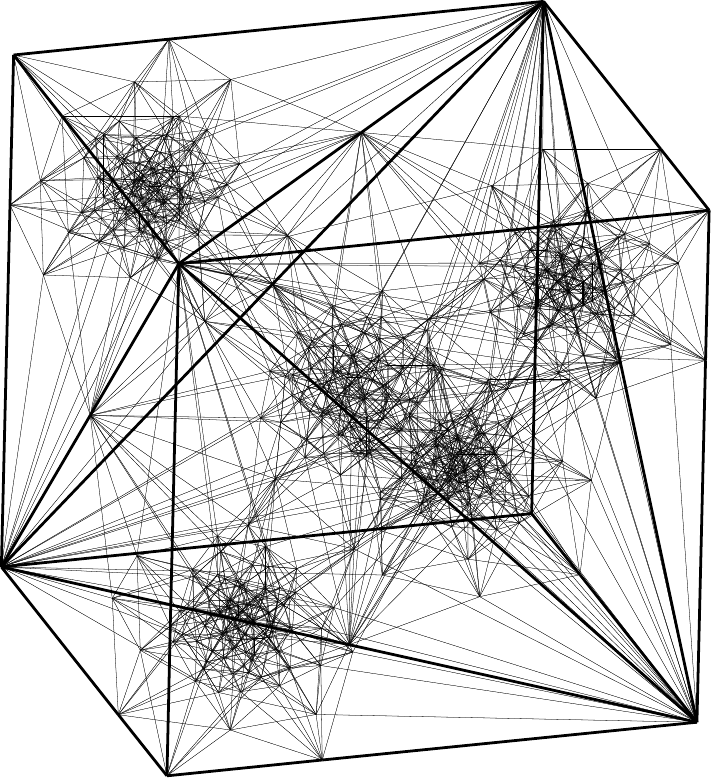}
\centerline{Fig. 7}
\end{center}

Finally, we obtain a triangulation satisfying conditions of
Theorem~\ref{dim3}, part~(2), by attaching eight copies of the
standard $3$--simplex triangulated as in Step~2. For further
discussion on the construction we refer to Subsection~\ref{final1}.

\begin{rem}
\label{Eryk}
The animation of our triangulation of the $3$--cube is available at \\
\centerline{\texttt{http://www.mimuw.edu.pl/\~{}erykk/papers/acute.html}.}
For the details and the exact values of all the parameters which
have been guessed, see the implementation of above algorithm,
available with the animation.
\end{rem}

\section{Dehn--Sommerville equations in dimension~$4$}
\label{section D-S}

In this short section we present some known results in geometric
combinatorics.

Denote by $f_i(M)$, $f_i(\partial M)$ (we later abbreviate this to
$f_i, f^\p_i$), the number of $i$--dimensional simplices of a
triangulation of a compact $m$--dimensional homology manifold $M$ 
and its boundary $\partial M$. Recall the following
Dehn--Sommerville type equations (see Subsection~\ref{DS} for the
history of this generalization).

\begin{theorem}[{\cite[Theorem 1.1]{Kl}} and \cite{NS}]
\label{D-S general} Let $M$ be a~compact~$m$--dimensional
triangulated homology manifold with boundary. For $k=0, \ldots, m$ 
we have
$$ f_k(M)-f_k(\partial M)=\sum_{i=k}^m(-1)^{i+m}{{i+1} \choose {k+1}}f_i(M).$$
\end{theorem}

If $m=4$, then for $k=1,2$ we obtain the following.

\begin{cor}
\label{D-S}
If $M$ is $4$--dimensional and we abbreviate $f_i=f_i(M), f^\p_i=f_i(\partial M)$, then
\item[(i)]$2f_1-f^\p_1=3f_2-6f_3+10f_4$,
\item[(ii)]$-f^\p_2=-4f_3+10f_4$.
\end{cor}

These equalities will be used repeatedly in the next two sections.

\section{Rich triangulations of $4$--manifolds}
\label{Rich triangulations of 4--manifolds}

In this section we prove the following combinatorial result on rich
triangulations (see Definition \ref{rich}) of $4$--dimensional
homology manifolds. This addresses Przytycki--\'Swi\k{a}tkowski 
\cite[Questions 5.8(3)]{PS}. We keep the notation $f_i, f^\p_i$ from
Section~\ref{section D-S}.

\begin{theorem}
\label{main} Every rich triangulation of a compact $4$--dimensional
homology manifold $M$ with Euler characteristic $\chi$ satisfies 
$$2f_0\leq 2\chi+f^\p_1.$$ In particular, if $M$ is closed, then $f_0\leq \chi$.
\end{theorem}

Before we present the proof of the theorem, let us give the
following four corollaries in the case when the homology
manifold~$M$ is closed. 

\begin{cor}
Any $4$--dimensional closed homology manifold~$M$ has only finitely 
many rich triangulations. In particular~$M$ has only finitely many
flag-no-square triangulations.
\end{cor}

\begin{cor}[Theorem \ref{2}]
\label{4-cube} There is no periodic (i.e.\ invariant under
a~cocompact group of translations) acute triangulation of $\R^4$. In
particular, there is no acute triangulation of the $4$--cube.
\end{cor}

\proof A~periodic triangulation $\tau$ of $\R^4$ descends to
a~triangulation $\tau'$ of a~$4$--torus. Since the Euler
characteristic of a $4$--torus equals $0$, by Theorem~\ref{main},
triangulation~$\tau'$ is not rich. Hence, by Observation \ref{five},
$\tau$ is not acute. This proves the first part of the corollary.
For the second part, observe that an acute triangulation of
the~$4$--cube could be promoted, by reflecting, to a~periodic acute
triangulation of $\R^4$. \qed

\medskip

Let us show also that the theorem gives a much simplified proof of
the following known result:

\begin{cor}[{\cite[Section 2.2]{JS}}]
There are no rich (in particular no flag-no-square) triangulations
of closed homology $n$--manifolds, for $n\geq 5$. 
\end{cor}

\proof The link $L$ of any codimension $5$ simplex of a
triangulation $\sigma$ of a closed homology $n$--manifold (for 
$n\geq 5$) is a $4$--dimensional homology sphere, which implies that
its Euler characteristic $\chi$ equals $2$. Since $L$ is
$4$--dimensional, it must have at least $6$ vertices. Hence, by
Theorem \ref{main}, $L$ is not rich. Thus, triangulation~$\sigma$ is
also not rich, which proves the result. \qed

\medskip

In the same way one we obtain the following theorem (cf. Subsection~\ref{Krizek}):

\begin{cor}[Theorem \ref{higher dim introduction}, {\cite[Theorem 6.2]{Kri}}]
\label{higher dim}
There is no triangulation of a polyhedron in $\R^n$, for $n\geq 5$, which
contains an interior vertex such that all dihedral angles adjacent to it
are acute.
\end{cor}
\proof
Let $v$ be an interior vertex and let $\rho$ be a codimension $5$ simplex (a vertex for $\R^5$, an edge for $\R^6$ etc) containing $v$.
The link $L$ of $\rho$ is a $4$--dimensional homology sphere,
hence its Euler characteristic equals $2$. Since $L$ is $4$--dimensional, it must have at least $6$ vertices.
Hence, by Theorem \ref{main}, $L$ is not rich. Thus the link of one of the codimension 2 simplices containing $\rho$
is a cycle of length shorter than $5$. Hence one of the dihedral angles adjacent to $v$ is not acute.
\qed

\medskip

Finally we provide the following.

\medskip\par\noindent\textbf{Proof of Theorem \ref{main}.}\ignorespaces
\ Let $\tau$ be a~rich triangulation of $M$. We compute the
number~$N$ of flags $(\rho_2 \subset \rho_4)$ of
a~$2$--simplex~$\rho_2$ contained in a~$4$--simplex~$\rho_4$. On one
hand, it equals $10\ts f_4$, since each $4$--simplex has ten
$2$--dimensional faces. On the other hand, by richness, each
interior $2$--simplex (there are $f_2-f^\p_2$ of those) is contained
in at least five $4$--simplices. Thus, we have:
\begin{equation}
\label{richeq}
N \. = \. 10\ts f_4 \. \geq \. 5 \ts (f_2-f^\p_2).
\end{equation}

By the definition of the Euler characteristic, we have:
\begin{equation}
\label{Euler}
\chi-f_0=-f_1+f_2-f_3+f_4.
\end{equation}
Applying consecutively formula~\eqref{Euler}, then Corollary~\ref{D-S} parts~(i) and~(ii),
and formula~\eqref{richeq}, we obtain:
\begin{align*}
2(\chi-f_0)+f^\p_1&=-2f_1+2(f_2-f_3+f_4)+f^\p_1=\\
                  &=-(f^\p_1+ 3f_2-6f_3+10f_4)+2(f_2-f_3+f_4)+f^\p_1=\\
                  &=-f_2+4f_3-8f_4=-f_2+(f^\p_2+10f_4)-8f_4=\\
                  &=2f_4-(f_2-f^\p_2) \geq 0,\\
\end{align*}
as desired.
\qed

\section{Acute triangulations of $\R^4$}
\label{Acute-angled triangulations of R4} In this section (see also
Section~\ref{remarks R4}) we address the problem whether there is an
acute triangulation of $\R^4$. We know already that every such acute
triangulation of $\R^4$ cannot be periodic (Corollary~\ref{4-cube}).
Here we present the following stronger result.

We say that a~triangulation of $\R^p$ has \emph{bounded geometry} if there is a global
upper~bound on the ratio of edge lengths in every $p$--simplex.

\begin{theorem}
\label{no in R^4}
There is no acute triangulation of $\R^4$ with bounded geometry.
\end{theorem}

This result can be restated in the following (equivalent) form:

\begin{cor}[Theorem \ref{3}]
\label{away from 90}
There is no acute triangulation of $\R^4$ with dihedral angles bounded away from $\frac{\pi}{2}$.
\end{cor}

\proof
If the dihedral angles are bounded away from $\frac{\pi}{2}$,
then the angles of $2$--simplices are bounded
away from $\frac{\pi}{2}$ (see e.g. \cite{Kri}).
Hence the angles of $2$--simplices are also bounded away from~$0$.
By the sine law, this gives a~bound on the ratio of lengths
of edges in each $2$--simplex, which results in a~bound of
the ratio of lengths of edges in each $4$--simplex.
\qed

\medskip

Before we prove Theorem \ref{no in R^4}, we need a few
preliminary results.

\begin{lemma}
\label{quasiconformal}
Let $\tau$ be a triangulation of $\R^p$ with bounded geometry.
Then the $1$--skeleton of~$\tau$ has bounded degree.
\end{lemma}
\proof
All $p$--simplices of $\tau$ are affinely quasi-conformal to the regular $p$--simplex with a~universal constant.
Hence the spherical volume contributed by any $p$--simplex in the link of any vertex of $\tau$
is bounded from below. On the other hand, the total volume of the link is the volume of the unit $(p-1)$--sphere.
This bounds the number of $p$--simplices, and in particular the number of lower dimensional simplices,
sharing each vertex.
\qed

\begin{lemma}
\label{estimating Euler} Let $M$ be a compact connected
$4$--dimensional triangulated homology manifold. Assume that $M$
admits a PL embedding into $\R^4$. Then the Euler characteristic of
$M$ is at most $1+\mathrm{rk} H_2(\partial M)$.
\end{lemma}
\proof
First observe that the natural map $H_2(M)\rightarrow H_2(M,\partial M)$ is trivial. Indeed, this mapping factors through
$$H_2(M)\rightarrow H_2(S^4)\rightarrow H_2(S^4, S^4\setminus M)=H_2(M,\partial M),$$ where $S^4$ is the one point
compactification of $\R^4$ with $H_2(S^4)=0$.
Hence the natural map $H_2(\partial M)\rightarrow H_2(M)$ is onto and
$\mathrm{rk} H_2(M)\leq \mathrm{rk} H_2(\partial M)$.
Thus, the Euler characteristic~$\chi$ of $M$ satisfies
$$
\chi \. \le \. \mathrm{rk}H_0(M)+\mathrm{rk} H_2(M)=1+\mathrm{rk} H_2(M)
\. \leq \. 1+\mathrm{rk} H_2(\partial M),
$$ as desired.
\qed
\medskip
Theorem~\ref{main} and Lemma~\ref{estimating Euler} now imply the
following result:

\begin{cor}
\label{comb corollary} Let $M$ be a compact connected
$4$--dimensional homology manifold with a rich triangulation. Assume
that $M$ admits a PL embedding into $\R^4$. If $f_i, f^\p_i$ are
defined as in Section~\ref{section D-S}, then we have
$$2f_0\leq 2(1 + f^\p_2) +f^\p_1.$$
\end{cor}

We turn our attention now to the study of isoperimetric functions on infinite graphs.

\begin{defin}
Let $G=(V,E)$ be a simple connected (locally finite) infinite~graph,
and let $\Omega\subset V$ be a finite subset of vertices. Denote by
$\partial \Omega$ the \emph{vertex-boundary} of~$\Omega$, defined as
the subset of $V\setminus \Omega$ consisting of vertices adjacent to
vertices in~$\Omega$.

We say that $I\colon \mathbb{Z}_{\geq 0}\rightarrow \mathbb{Z}_{\geq
0}$ is an \emph{isoperimetric function} for~$G$, if the inequality
$I(|\Omega|)\leq |\partial \Omega|$ holds for every finite
$\Omega\subset V$.
\end{defin}

\begin{proposition}
\label{isoperimetry}
The $1$--skeleton of any acute triangulation of $\R^4$ with bounded geometry has linear isoperimetric function.
\end{proposition}
\proof Let $G=(V,E)$ be the $1$--skeleton of an acute triangulation
of $\R^4$ with bounded geometry. Consider any finite $\Omega\subset
V$. We want to obtain a linear isoperimetric function for $G$, hence
we may assume that the subgraph spanned by $\Omega$ in $V$ is
connected. 

To outline the idea of the proof, assume first that the subcomplex
of $\R^4$ which is the closure of the union of all simplices meeting
$\Omega$ is a~($4$--dimensional) homology manifold. Denote this
subcomplex by $M$. Then $\Omega$ is contained in $M$ and the
vertices in $\partial M$ lie in $\partial \Omega$. By Lemma
\ref{quasiconformal}, $f^\p_1$ and $f^\p_2$ are bounded above by
$Cf^\p_0$, for some fixed $C$. Hence, by Corollary \ref{comb
corollary}, we have
$$ |\Omega|\leq f_0 \leq 1+f^\p_2+\frac{1}{2}f^\p_1\leq 1+\frac{3}{2}Cf^\p_0
\leq (1+\frac{3}{2}C)f^\p_0\leq (1+\frac{3}{2}C)|\p \Omega|,$$ as
desired.

In general, as pointed out to us by Jon McCammond, the closure of
the union of all simplices meeting $\Omega$ might not be a homology
manifold. The strategy then is, roughly speaking, to subdivide the
original triangulation in order to find a tubular PL neighborhood
$M$ of $\Omega$, which is a homology $4$--manifold with rich
triangulation. To this end, we need the following construction:

\begin{defin}
Let $Y$ be a subcomplex of a simplicial complex $X$. Let $N_X(Y)$ be
the simplicial complex containing $Y$ defined in the following way.
Its set of vertices is the union of the set of vertices of $Y$ and
of the set of simplices of $X$ which are not contained in $Y$, but
contain a vertex from $Y$. Now we describe the simplices in
$N_X(Y)$. Assume that a simplex $\sigma$ is contained in $Y$ and we
have simplices $\sigma\subset \tau_1\subset\ldots \subset \tau_k$ in
$X$ with $\tau_1 \nsubseteq Y$. Then in $N_X(Y)$ we span a simplex
on the union of the set of vertices from $Y$ lying in $\sigma$ and
on the set $\{\tau_i\}$.
\end{defin}

We assume now that $X$ is a homology $n$--manifold without boundary.
Suppose that $Y$ is a \emph{full} subcomplex of $X$ (i.e.\ if all
the vertices of a simplex $\sigma$ of $X$ belong to $Y$, then also
$\sigma$ belongs to $Y$). Then the simplicial complex $M=N_X(Y)$ is
also a homology $n$--manifold. The image of the natural embedding of
$M=N_X(Y)$ into $X$, which restricts to the identity on $Y$ and maps
each vertex corresponding to a simplex of $X$ to its barycenter in
$X$, can be regarded as a PL tubular neighborhood of $Y$ in $X$.

The vertices in the interior of $M$ are exactly the vertices of $Y$.
The link in $M$ of a $(n-2)$--dimensional simplex in $Y$ is obtained
by subdividing its link in $X$. Hence if $X$ is rich, then $M$ is
also rich.

\medskip
We now return to the proof of Proposition~\ref{isoperimetry}. Let
$X$ be the $\R^4$ equipped with the acute triangulation with bounded
geometry and let $Y\subset X$ be the full subcomplex spanned by
$\Omega$. Put $M=N_X(Y)$. Since $Y$ is connected, $M$ is connected
as well. Since $\Omega$ is the set of vertices of $Y$, we have
$|\Omega|\leq f_0$. On the other hand, every vertex in the boundary
of $M$ corresponds to a simplex in $X$ containing a vertex of
$\partial \Omega$. Hence by Lemma \ref{quasiconformal} we have
$f^\p_0\leq
D|\p \Omega|$, for some fixed $D$. 

Moreover, again by Lemma \ref{quasiconformal}, $f^\p_1$ and $f^\p_2$
are bounded above by $Cf^\p_0$, for some fixed $C$. Finally, by
Corollary \ref{comb corollary}, we have
$$ |\Omega|\leq f_0 \leq 1+f^\p_2+\frac{1}{2}f^\p_1\leq 1+\frac{3}{2}Cf^\p_0
\leq (1+\frac{3}{2}C)f^\p_0\leq (1+\frac{3}{2}C)D|\p \Omega|,$$ as
desired. \qed

\medskip
Summarizing, we showed that acute triangulations with bounded
geometry have a linear isoperimetric function. In the remaining part
of this section, we show that this leads to a~contradiction. The
argument that follows was suggested to us by Marc Bourdon. Following
Benjamini--Curien \cite[Section 2.2]{BC}, we recall the following
definition:

\begin{defin}
Let $G=(V,E)$ be a~locally finite connected graph and let $\Gamma(v)$ be the set of all semi-infinite self
avoiding simplicial paths in $G$
starting from $v\in V$. For any $m\colon V\rightarrow \R_+$ (so called \emph{metric}), the \emph{length} of a~path
$\gamma$ in $G$ is defined by
$$\text{Length}_m(\gamma)=\sum_{v\in \gamma}m(v).$$

If $m\in L^p(V)$, we denote by $||m||_p$ the usual $L^p$ norm. The graph $G$ is \emph{$p$--parabolic} if
the \emph{$p$--extremal length} of $\Gamma(v)$,
$$\sup_{m \in L^p(V)} \inf _{\gamma \in \Gamma(v)} \frac{\text{Length}_m(\gamma)^p}{||m||_p^p}$$
is infinite.
This definition does not depend on the choice of $v\in V$.
\end{defin}

\begin{lemma}
\label{R parabolic}
Let $G$ be the $1$--skeleton of a triangulation of $\R^p$ with bounded geometry, where $p\geq 2$.
Then $G$ is $p$--parabolic.
\end{lemma}

This lemma can be obtained from the Bonk and Kleiner
result~\cite[Corollary 8.8]{BK}. To make the proof complete and
self-contained, we include a concise proof.

\proof
Let $\ell,L\colon V\rightarrow \R_+$ be the length functions of the shortest and the longest edge
adjacent to a~vertex.
Since our triangulation has bounded geometry, by Lemma~\ref{quasiconformal} there is
a~constant $C>0$, such that $L(v) \leq C \ts \ell(v)$ for all $v \in V$.

We fix a~basepoint vertex $v\in V$. Let $m\colon V \rightarrow \R_{\geq 0}$ be a function
defined by
$$m(v) \, = \, \frac{\ell(w)}{||w-v||} \qquad \text{for all} \ \, w \in V, \, w \ne v\ts,$$
and let $m(v)=0$.
For every $R\geq 0$, define
$m_R\colon V \rightarrow \R_{\geq 0}$ by $m_R(w)=m(w)$, for all $w \in V\cap B_{e^R}(v)$,
and $m_R(w)=0$ for all $w\notin B_{e^R}(v)$. Here and throughout this section,
$B_t(v)$ denotes the closed ball in~$G$ of distance~$t$, around vertex~$v$.
We claim that
\begin{equation}
\label{extremal} \inf _{\gamma \in \Gamma(v)}
\frac{\text{Length}_{m_R}(\gamma)^p}{||m_R||_p^p}\rightarrow \infty
\qquad \text{as} \ \. R\rightarrow \infty \ts,
\end{equation}
and therefore, the $p$--extremal length of $\Gamma(v)$ is infinite.

\medskip\par\noindent\textbf{Step 1.}\ignorespaces
\ Let $\gamma\in \Gamma(v)$. We begin with bounding
$\text{Length}_{m_R}(\gamma)$ from below. Let $\omega$ be the
$1$--form on $\R^p$ which is zero on $B_{L(v)}(v)$ and equal
$\frac{dr}{r}$ outside $B_{L(v)}(v)$, where $r$ is the radial
coordinate w.r.t. the basepoint $v$. Let $\gamma_R$ be the maximal
initial part of $\gamma$ consisting of vertices in $B_{e^R}(v)$ and
edges starting at vertices in $B_{e^R}(v)$. Since $\gamma_R$ starts
at $v$ and eventually leaves $B_{e^R}(v)$, we have that
$\int_{\gamma_R}\omega \geq R - \ln L(v)$. On the other hand, for an
edge $f$ of $\gamma_R$ starting at $w\neq v$ we have
$$\int_{f}\omega \leq \frac{L(w)}{||w-v||}.$$ Altogether, we obtain:

\begin{align*}
\text{Length}_{m_R}(\gamma) \, = \, & \sum _{w\in \gamma} \. m_R(w) \, \geq \,
\sum _{w\in \gamma_R}\. m_R(w) \, = \, \sum _{w\in \gamma_R} \. m(w)\, \geq\\
\geq \, & \frac{1}{C}\sum _{w\in \gamma_R\setminus v}\.
\frac{L(w)}{||w-v||} \, \geq \frac{1}{C}\. \sum_{f\in \gamma_R} \.
\int_f \omega \, =
\, \frac{1}{C}\int_{\gamma_R}\omega \, \geq \, \frac{1}{C}\ts (R - \ln L(v)).\\
\end{align*}

\medskip
\par\noindent\textbf{Step 2.}\ignorespaces
\ We now bound $||m_R||_p$ from above. With each vertex $w\in V$ we associate the ball $B(w)$ of radius
$\frac{l(w)}{2}$ centered at $w$. All these balls have disjoint interiors.

Let $\sigma$ be the $p$--form on $\R^p$ which is zero on $B(v)$ and is equal to
$\frac{1}{r^p}\text{vol}$ outside~$B(v)$, where \ts vol \ts
denotes the Euclidean volume form.

Let $w\neq v$. We estimate $\sigma(B(w))$. Since the radius of $B(w)$ is at most $\frac{||w-v||}{2}$, we have
$B(w)\subset B_{\frac{3||w||}{2}}(v)$. Hence:
\begin{align*}
\sigma(B(w))=& \int_{B(w)}\frac{\text{vol}}{r^p}\geq \left(\frac{2}{3||w-v||}\right)^p\text{vol}(B(w))=\\
=&\left(\frac{2}{3||w-v||}\right)^p V_p \left(\frac{l(w)}{2}\right)^p=c\ts m^p(w),\\
\end{align*}
for a~universal constant $c>0$, where $V_p$ is the volume
of the unit $p$--ball. Hence
$$ ||m_R||^p_p =\sum _{w\in B_{e^R}(v)}m^p(w)\leq \frac{1}{c} \sum_{w\in B_{e^R}(v)}\sigma(B(w))\leq
\frac{1}{c}\sigma(B_{\frac{3}{2}e^R}(v)).$$ The latter is bounded above by
$A_{p-1}(\ln \frac{3}{2}+R-\ln \frac{l(v)}{2})$,
where $A_{p-1}$ is the volume
of the unit $(p-1)$--sphere,
and that is a~linear function in $R$.

\medskip\par\noindent\textbf{Combining the steps.}\ignorespaces
\ In Step 1 we have bounded the numerator of \eqref{extremal} below by a~polynomial of degree $p$ in $R$. In
Step 2 we have bounded the denominator of \eqref{extremal} above by a~function linear in $R$. Hence
the $p$--extremal length of $\Gamma(v)$ is infinite, as desired. \qed

\bigskip
To finish the proof, we need the following known result:

\begin{proposition}[{\cite[Proposition 4.1(1)]{BC}}]
\label{Maillot}
Let $G=(V,E)$ be an infinite locally finite connected graph. If $G$ is $p$--parabolic and $I$ is an
isoperimetric function, then for $D>p$ we have
$$\sum^\infty_{k=1}\frac{1}{I(k)^\frac{p}{p-1}}=\infty.$$
\end{proposition}

\medskip\par\noindent\textbf{Proof of Theorem \ref{no in R^4}.}\ignorespaces
\ Assume that there is is an acute triangulation $\tau$ of $\R^4$
with bounded geometry. Let $G$ be the $1$--skeleton of $\tau$. By
Lemma \ref{R parabolic} we have that the graph~$G$ is
$4$--parabolic. By Proposition \ref{Maillot}, we have that
$k\rightarrow Ck$ is not an isoperimetric function for~$G$. This
contradicts Proposition~\ref{isoperimetry}. \qed

\section{Final remarks and open problems}
\label{final}

\subsection{}
\label{final1}
It is unclear how far the results of Section~\ref{construction} extend
to other polytopes in~$\R^3$. For example, the regular icosahedron has an easy acute
triangulation using cones from the center over every facet. Similarly, the regular
dodecahedron, one can easily subdivide it into 120 congruent tetrahedra all meeting
at the center. It turns out that the special subdivision of this triangulation has
an acute realization in this case as well. For the subdivision of one of the
congruent tetrahedra see\\
\centerline{\texttt{http://www.mimuw.edu.pl/\~{}erykk/papers/acute.html}.}
Putting this together, we obtain the following result:

\begin{theorem}
All Platonic solids have non-trivial acute triangulations.
\end{theorem}

Now, it is possible that \emph{every} convex polytope in~$\R^3$ has an acute
triangulation. We conjecture this to be the case. Unfortunately, we are very
far from proving this, given that this paper and~\cite{VH} have the first examples
of non-trivial acute triangulations of convex polytopes (cf.~\cite{BS}).
Perhaps, it is even possible that every $3$--dimensional abstract polyhedral manifold
has an acute triangulation,
in the style of~\cite{BZ}. For example, we conjecture that the boundary of every
convex polytope in~$\R^4$ has an acute triangulation. Of course, in the spirit
of~\cite{BGR,M,Sar}, the problem might prove much easier for non-obtuse triangulation.

Finally, numerical results would also be of interest. What is the smallest number
of tetrahedra required for a non-trivial acute triangulation of the regular
tetrahedron? For example, can one beat our record of~543? How about the cube?
Can one bound the smallest maximum dihedral angle? Dreaming of the future, can one
always find a linear size acute triangulation of a convex polytope in~$\R^3$?

\subsection{}
\label{remarks R4}
Although K{\v{r}}{\'{\i}}{\v{z}}ek conjectured in~\cite{Kri} (see
also~\cite{BS}), that there are no
acute triangulations of the space~$\R^4$, our results resolve only a special case
of this problem. The conjecture remains open in full generality, when the geometry
of simplices is unbounded. On the one hand, another (plausible)
conjecture in~\cite{Kri,BS} states that locally such acute triangulation
must have at least $600$ simplices around each vertex, making a construction of such
triangulation exceedingly difficult. On the other hand, in the plane and the space
there are known very general combinatorial tiling constructions which require an unbounded
geometry (see e.g.~\cite{GMS,Sch}). We conjecture that there exists an acute triangulation
of~$\R^4$, although we think that to construct it, one first has to master acute
triangulations in dimension $3$ (see Section~\ref{final1}), and in the spherical
case (see Section~\ref{spherical} below).

Similarly, what happens with individual convex polytopes in~$\R^4$ is much less clear,
and will obviously depend on the polytope. For example, by analogy with the icosahedron,
there is an easy acute triangulation of the $600$--cell. On the other hand, it is unclear
whether the $16$--cell (the regular cross-polytope),
the $24$--cell and the $120$--cell have acute triangulations (we conjecture not).
One is tempted to conclude the $16$--cell does not admit an acute triangulation, since it
tiles the space~$\R^4$. Unfortunately, this argument is incorrect for the following reason.
In order to have consistency on the boundary, two subdivisions of a $16$--cell adjacent by
the tetrahedra must have the opposite orientations. However, in the tiling, there are three
(an odd number) top dimensional cells around each codimension $2$ simplex, and thus not
every triangulation of the $16$--cell gives rise to a triangulation of $\R^4$.

Interestingly, the space tessellation argument does work for some
notable polytopes in~$\R^4$. For example, it is well known that the
$4$--cube can be dissected into $24$ congruent
orthoschemes\footnote{They are also called path-simplices.} (see
e.g.~\cite{C}), in such a way that around each interior codimension
$2$ simplex there are~$4$ or~$6$ (an even number) orthoschemes. This
implies that such ``isosceles'' orthoscheme does not have an acute
triangulation, because otherwise one could extend it by reflecting
to an acute triangulation of the $4$--cube. Similarly, the standard
(cube corner) $4$--dimensional simplex tiles the space $\R^4$ by
reflections (it is a fundamental chamber for $\widetilde D_4$ affine
Coxeter group action). Hence the standard
$4$--simplex also does not admit an acute triangulation. 

In summary, we believe that finding a useful criterion for a polytope
in~$\R^4$ to have an acute triangulation is a challenging problem, which
we expect to be more difficult than the $3$--dimensional version.

\subsection{}
\label{Krizek}
We believe that the proof of our Theorem~C given in~\cite[Theorem~6.2]{Kri}
has a crucial gap and is either incomplete or incorrect as
written.\footnote{After this paper was written, Michal
K{\v{r}}{\'{\i}}{\v{z}}ek graciously confirmed the error in his paper.
He informed us that he first learned about it in 2008 from Jan Brandts,
and that he has prepared a correction (to appear). Since neither the
error nor the correction has been announced nor are publicly available,
we decided not to change our presentation and keep the details.}
On p.~387, in the proof of Theorem 5.1, in the sentence \ts ``the sum
of all dihedral angles of tetrahedra around a~given edge~$E$
from $\partial P$ cannot be greater than~$2\pi$,'' \ts the author seems to be
referring to $3$--faces of convex $4$--polytope~$P$, which would make this
statement true (and Lemma~3.3 in the paper applicable). However, throughout
the paper, the polytope~$P$ is in fact in~$\R^5$, in which case the above
sum is a priori unbounded. It seems, this mistake has not been discovered until now.
We should mention, however, that the reduction of higher dimensions to dimension~5
given in~\cite{Kri} (see the proof of Theorem~6.2), is independent of Theorem~5.1
and completely correct.

\subsection{}
\label{spherical}
It would be interesting to consider the spherical and hyperbolic analogues
of the acute triangulation problem. The spherical analogue might prove particularly
insightful as it might allow the use a dimension reduction in the Euclidean case
(in particular in the case of~$\R^4$).

\subsection{} The Burago--Zalgaller original result in~\cite{BZ} is a technical lemma
used towards the $3$--dimensional analogue of the classical Nash--Kuiper embedding theorem.
This result (by a different technique) was recently extended by Akopyan to higher
dimensions~\cite{Ako}, despite the absence of acute triangulations.

\subsection{} In the plane, one can start with a given triangulation and ``improve it''
by using $2$--flips, by increasing the smallest angle in a triangle.
This results in the \emph{Delaunay triangulation} which (among all
triangulations on this set of vertices) has the largest possible
minimal angle, and has a number of other useful properties (though
not necessarily the smallest maximal angle). Thus, by strategically
placing new points into the interior of a polygon one can then
efficiently construct a ``good'' triangulation. In higher dimensions
this approach breaks down for several reasons, both due to the lack
of connectedness of triangulations via flips, and non-monotonicity
of the angle functionals. We refer to~\cite{DRS,Pak} for an
introduction and an extended discussion of the problem.

Interestingly, a variation on the Delaunay approach does give useful meshes in~$\R^3$,
as described in~\cite{VR}. The paper~\cite{VH} is a followup on this approach, which
uses a more refined idea of incremental changes in a triangulation, by moving the
vertices one at a time.

\subsection{} There is a large body of work on tessellations of the space by
convex polyhedra with bounded geometry. Perhaps the earliest, is the
result of Alexandrov that in every triangulation of the plane into
bounded triangles the average degree of vertices (when defined) must
be at least~6~\cite{Alex}. Another is a classical result by Niven
that convex $n$--gons of bounded geometry cannot tile the plane
for~$n \ge 7$~\cite{Niven} (see also~\cite{Fulton}). The idea
is always to use the isoperimetric inequalities compared with direct
counting estimates, an approach which works in higher dimensions
as well (see e.g.~\cite{BSch,KS,LM}). Our proof of Theorem~D is a
variation on the same line of argument. We refer to~\cite{GS} for
historical background and further references.

\subsection{}
\label{DS}
The classical Dehn--Sommerville equations are defined for $f$--vectors of
simplicial convex polytopes in~$\R^d$ and (see e.g.~\cite[Section~8]{Pak}).
They are extended in a number of directions, notably the beautiful
\emph{flag $f$--vectors} by Bayer and Billera, leading to the \emph{$cd$-index}
(see~\cite{St-cd}). The version for manifolds without boundary was first given
by Klee in~\cite{Klee}. It seems, the version with the boundary goes back to
Macdonald~\cite{Mac} and was rediscovered a number of times. We refer
to~\cite{NS} for the most general version of the equations, various applications,
and further references.

\subsection{}
\label{exact}
Let us note that the ad hoc acute triangulation of the cube
discovered in~\cite{VH} has 1370 tetrahedra as opposed to $5\cdot 543 = 2715$
tetrahedra in our construction. On the other hand, one can argue that our
construction is more symmetric, including the natural action of~$S_4$ by permuting the
standard tetrahedra, as well as some ``hidden'' symmetries arising from the
$600$--cell (the construction in~\cite{VH} also has a number of symmetries).

After we made our computations publicly available, Evan VanderZee
kindly informed us that he re-checked our coordinates for the
triangulation of the cube and computed that dihedral angles range
between 26.425 and~89.992.\footnote{It should be noted that
optimizing the angles was not our goal. The reason the angles are so
large, is because our angle correction procedure described in
Section~\ref{construction} works only with angles larger or equal
than~$\pi/2$.} For comparison, the dihedral angles found
in~\cite{VH} range between 35.89 and~84.65, which is significantly
better for numerical algorithms, since it has fewer tetrahedra,
\emph{smaller} the maximal and \emph{larger} the minimal dihedral
angles. On the other hand, after performing simulations with our
mesh of 2715 tetrahedra (i.e.\ when combinatorial structure is fixed
while positions of points are allowed to vary), VanderZee obtained
an acute triangulation with dihedral angles between 25.310
and~88.902. Hence this triangulation and ours are incomparable as
far as the dihedral angles are concerned.

\appendix

\section{The exact position of vertices}

In order for the coordinates to be integers, we triangulate the cube whose vertices are at the eight points whose each
coordinate is 0 or 60000 (instead of $\pm 1$). Vertices of the Step~3 (Section~\ref{construction}) triangulation
of the regular $3$--simplex $T_0$ have the following coordinates (we list four vertices in each row):

\def\posarray#1{{
  \tiny
  \begin{longtable}{@{\extracolsep{-.25em}}lrrrrrrrrrrrr}
  #1
  \end{longtable}
  }}

\def\line#1#2#3#4#5{#1 & #2; & #3; & #4; & #5 \\}
\def\pt#1#2#3{#1,&#2,&#3}

\posarray{
\line{0-3}{\pt{60000}{0}{0}}{\pt{60000}{60000}{60000}}{\pt{0}{0}{60000}}{\pt{0}{60000}{0}}
\line{4-7}{\pt{0}{30000}{30000}}{\pt{60000}{30000}{30000}}{\pt{30000}{30000}{0}}{\pt{30000}{60000}{30000}}
\line{8-11}{\pt{30000}{30000}{60000}}{\pt{30000}{0}{30000}}{\pt{33916}{43042}{16958}}{\pt{43042}{26084}{43042}}
\line{12-15}{\pt{16958}{16958}{26084}}{\pt{16958}{33916}{43042}}{\pt{43042}{16958}{33916}}{\pt{43042}{33916}{16958}}
\line{16-19}{\pt{43042}{43042}{26084}}{\pt{16958}{26084}{16958}}{\pt{16958}{43042}{33916}}{\pt{26084}{43042}{43042}}
\line{20-23}{\pt{26084}{16958}{16958}}{\pt{33916}{16958}{43042}}{\pt{34171}{39326}{34171}}{\pt{20674}{34171}{25829}}
\line{24-27}{\pt{25829}{25829}{39326}}{\pt{34171}{25829}{20674}}{\pt{25829}{20674}{34171}}{\pt{39326}{34171}{34171}}
\line{28-31}{\pt{39326}{25829}{25829}}{\pt{20674}{25829}{34171}}{\pt{34171}{20674}{25829}}{\pt{34171}{34171}{39326}}
\line{32-35}{\pt{25829}{34171}{20674}}{\pt{25829}{39326}{25829}}{\pt{24956}{35044}{35044}}{\pt{39033}{32132}{27868}}
\line{36-39}{\pt{27868}{27868}{20967}}{\pt{32132}{20967}{32132}}{\pt{32132}{32132}{20967}}{\pt{24956}{24956}{24956}}
\line{40-43}{\pt{32132}{39033}{27868}}{\pt{20967}{27868}{27868}}{\pt{39033}{27868}{32132}}{\pt{35044}{35044}{24956}}
\line{44-47}{\pt{20967}{32132}{32132}}{\pt{27868}{20967}{27868}}{\pt{35044}{24956}{35044}}{\pt{27868}{32132}{39033}}
\line{48-51}{\pt{32132}{27868}{39033}}{\pt{27868}{39033}{32132}}{\pt{35393}{30000}{35393}}{\pt{24607}{30000}{24607}}
\line{52-55}{\pt{35393}{24607}{30000}}{\pt{24607}{35393}{30000}}{\pt{35393}{30000}{24607}}{\pt{24607}{30000}{35393}}
\line{56-59}{\pt{24607}{24607}{30000}}{\pt{35393}{35393}{30000}}{\pt{30000}{35393}{24607}}{\pt{30000}{24607}{35393}}
\line{60-63}{\pt{30000}{24607}{24607}}{\pt{30000}{35393}{35393}}{\pt{33844}{26156}{26156}}{\pt{33844}{33844}{33844}}
\line{64-67}{\pt{26156}{26156}{33844}}{\pt{26156}{33844}{26156}}{\pt{24207}{28632}{31368}}{\pt{31368}{31368}{35793}}
\line{68-71}{\pt{28632}{35793}{28632}}{\pt{28632}{28632}{35793}}{\pt{28632}{24207}{31368}}{\pt{35793}{31368}{31368}}
\line{72-75}{\pt{24207}{31368}{28632}}{\pt{35793}{28632}{28632}}{\pt{31368}{35793}{31368}}{\pt{31368}{28632}{24207}}
\line{76-79}{\pt{28632}{31368}{24207}}{\pt{31368}{24207}{28632}}{\pt{31992}{31992}{25546}}{\pt{34454}{28008}{31992}}
\line{80-83}{\pt{28008}{25546}{28008}}{\pt{28008}{31992}{34454}}{\pt{31992}{34454}{28008}}{\pt{25546}{28008}{28008}}
\line{84-87}{\pt{25546}{31992}{31992}}{\pt{34454}{31992}{28008}}{\pt{28008}{34454}{31992}}{\pt{28008}{28008}{25546}}
\line{88-91}{\pt{31992}{28008}{34454}}{\pt{31992}{25546}{31992}}{\pt{32775}{32775}{30655}}{\pt{27225}{30655}{27225}}
\line{92-95}{\pt{29345}{27225}{32775}}{\pt{32775}{29345}{27225}}{\pt{27225}{32775}{29345}}{\pt{27225}{29345}{32775}}
\line{96-99}{\pt{27225}{27225}{30655}}{\pt{32775}{30655}{32775}}{\pt{32775}{27225}{29345}}{\pt{30655}{27225}{27225}}
\line{100-103}{\pt{30655}{32775}{32775}}{\pt{29345}{32775}{27225}}{\pt{28494}{31506}{31506}}{\pt{30865}{29135}{29135}}
\line{104-107}{\pt{33159}{30000}{30000}}{\pt{28494}{28494}{28494}}{\pt{29135}{29135}{30865}}{\pt{26841}{30000}{30000}}
\line{108-111}{\pt{30000}{30000}{33159}}{\pt{30000}{26841}{30000}}{\pt{31506}{31506}{28494}}{\pt{29135}{30865}{29135}}
\line{112-115}{\pt{30000}{30000}{26841}}{\pt{31506}{28494}{31506}}{\pt{30865}{30865}{30865}}{\pt{30000}{33159}{30000}}
}

\noindent Vertices of the triangulation of the standard $3$--simplex
$T_1$ defined in Step 2 of Section~\ref{construction}, have the
following coordinates (rotate to obtain the coordinates for the
other standard~$T_i$):

\posarray{
\line{0-3}{\pt{60000}{0}{0}}{\pt{0}{0}{0}}{\pt{0}{0}{60000}}{\pt{0}{60000}{0}}
\line{4-7}{\pt{0}{30000}{30000}}{\pt{17574}{0}{0}}{\pt{30000}{30000}{0}}{\pt{0}{17574}{0}}
\line{8-11}{\pt{0}{0}{17574}}{\pt{30000}{0}{30000}}{\pt{12384}{20726}{0}}{\pt{10445}{0}{10445}}
\line{12-15}{\pt{16958}{16958}{26084}}{\pt{0}{12384}{20726}}{\pt{20726}{0}{12384}}{\pt{20726}{12384}{0}}
\line{16-19}{\pt{10445}{10445}{0}}{\pt{16958}{26084}{16958}}{\pt{0}{20726}{12384}}{\pt{0}{10445}{10445}}
\line{20-23}{\pt{26084}{16958}{16958}}{\pt{12384}{0}{20726}}{\pt{6257}{10104}{6257}}{\pt{9498}{21496}{13569}}
\line{24-27}{\pt{7743}{7743}{19656}}{\pt{21496}{13569}{9498}}{\pt{13569}{9498}{21496}}{\pt{10104}{6257}{6257}}
\line{28-31}{\pt{19656}{7743}{7743}}{\pt{9498}{13569}{21496}}{\pt{21496}{9498}{13569}}{\pt{6257}{6257}{10104}}
\line{32-35}{\pt{13569}{21496}{9498}}{\pt{7743}{19656}{7743}}{\pt{5137}{13344}{13344}}{\pt{15349}{8879}{5284}}
\line{36-39}{\pt{17685}{17685}{11260}}{\pt{16563}{7777}{16563}}{\pt{16563}{16563}{7777}}{\pt{16026}{16026}{16026}}
\line{40-43}{\pt{8879}{15349}{5284}}{\pt{11260}{17685}{17685}}{\pt{15349}{5284}{8879}}{\pt{13344}{13344}{5137}}
\line{44-47}{\pt{7777}{16563}{16563}}{\pt{17685}{11260}{17685}}{\pt{13344}{5137}{13344}}{\pt{5284}{8879}{15349}}
\line{48-51}{\pt{8879}{5284}{15349}}{\pt{5284}{15349}{8879}}{\pt{10649}{6407}{10649}}{\pt{13477}{17719}{13477}}
\line{52-55}{\pt{16305}{7821}{12063}}{\pt{7821}{16305}{12063}}{\pt{16305}{12063}{7821}}{\pt{7821}{12063}{16305}}
\line{56-59}{\pt{13477}{13477}{17719}}{\pt{10649}{10649}{6407}}{\pt{12063}{16305}{7821}}{\pt{12063}{7821}{16305}}
\line{60-63}{\pt{17719}{13477}{13477}}{\pt{6407}{10649}{10649}}{\pt{17102}{11055}{11055}}{\pt{9039}{9039}{9039}}
\line{64-67}{\pt{11055}{11055}{17102}}{\pt{11055}{17102}{11055}}{\pt{10544}{14025}{16177}}{\pt{8666}{8666}{12147}}
\line{68-71}{\pt{9383}{15016}{9383}}{\pt{9383}{9383}{15016}}{\pt{14025}{10544}{16177}}{\pt{12147}{8666}{8666}}
\line{72-75}{\pt{10544}{16177}{14025}}{\pt{15016}{9383}{9383}}{\pt{8666}{12147}{8666}}{\pt{16177}{14025}{10544}}
\line{76-79}{\pt{14025}{16177}{10544}}{\pt{16177}{10544}{14025}}{\pt{13876}{13876}{8805}}{\pt{13231}{8160}{11294}}
\line{80-83}{\pt{14921}{12984}{14921}}{\pt{8160}{11294}{13231}}{\pt{11294}{13231}{8160}}{\pt{12984}{14921}{14921}}
\line{84-87}{\pt{8805}{13876}{13876}}{\pt{13231}{11294}{8160}}{\pt{8160}{13231}{11294}}{\pt{14921}{14921}{12984}}
\line{88-91}{\pt{11294}{8160}{13231}}{\pt{13876}{8805}{13876}}{\pt{10992}{10992}{9324}}{\pt{12447}{15145}{12447}}
\line{92-95}{\pt{11891}{10223}{14589}}{\pt{14589}{11891}{10223}}{\pt{10223}{14589}{11891}}{\pt{10223}{11891}{14589}}
\line{96-99}{\pt{12447}{12447}{15145}}{\pt{10992}{9324}{10992}}{\pt{14589}{10223}{11891}}{\pt{15145}{12447}{12447}}
\line{100-103}{\pt{9324}{10992}{10992}}{\pt{11891}{14589}{10223}}{\pt{10088}{12458}{12458}}{\pt{13197}{11836}{11836}}
\line{104-107}{\pt{12891}{10406}{10406}}{\pt{13247}{13247}{13247}}{\pt{11836}{11836}{13197}}{\pt{11234}{13719}{13719}}
\line{108-111}{\pt{10406}{10406}{12891}}{\pt{13719}{11234}{13719}}{\pt{12458}{12458}{10088}}{\pt{11836}{13197}{11836}}
\line{112-115}{\pt{13719}{13719}{11234}}{\pt{12458}{10088}{12458}}{\pt{11382}{11382}{11382}}{\pt{10406}{12891}{10406}}
}

\medskip

\noindent
In both cases, the edges are spanned on the following pairs of vertices:
\bigskip

\noindent
{
\tiny
\begin{spacing}{0.5}
4-2, 4-3, 5-0, 5-1, 6-0, 6-3, 7-1, 7-3, 8-1, 8-2, 9-0, 9-2, 10-3, 10-6, 10-7, 11-1, 11-5, 11-8, 12-2, 12-4, 12-9, 13-2, 13-4, 13-8, 14-0, 14-5, 14-9, 14-11, 15-0, 15-5, 15-6, 15-10, 16-1, 16-5, 16-7, 16-10, 16-15, 17-3, 17-4, 17-6, 17-12, 18-3, 18-4, 18-7, 18-13, 19-1, 19-7, 19-8, 19-13, 19-18, 20-0, 20-6, 20-9, 20-12, 20-17, 21-2, 21-8, 21-9, 21-11, 21-14, 22-1, 22-7, 22-16, 22-19, 23-3, 23-4, 23-17, 23-18, 24-2, 24-8, 24-13, 24-21, 25-0, 25-6, 25-15, 25-20, 26-2, 26-9, 26-12, 26-21, 26-24, 27-1, 27-5, 27-11, 27-16, 27-22, 28-0, 28-5, 28-14, 28-15, 28-25, 29-2, 29-4, 29-12, 29-13, 29-24, 29-26, 30-0, 30-9, 30-14, 30-20, 30-25, 30-28, 31-1, 31-8, 31-11, 31-19, 31-22, 31-27, 32-3, 32-6, 32-10, 32-17, 32-23, 33-3, 33-7, 33-10, 33-18, 33-23, 33-32, 34-13, 34-18, 34-19, 35-5, 35-15, 35-16, 35-27, 35-28, 36-6, 36-17, 36-20, 36-25, 36-32, 37-9, 37-14, 37-21, 37-26, 37-30, 38-6, 38-10, 38-15, 38-25, 38-32, 38-36, 39-12, 39-17, 39-20, 39-36, 40-7, 40-10, 40-16, 40-22, 40-33, 41-4, 41-12, 41-17, 41-23, 41-29, 41-39, 42-5, 42-11, 42-14, 42-27, 42-28, 42-35, 43-10, 43-15, 43-16, 43-35, 43-38, 43-40, 44-4, 44-13, 44-18, 44-23, 44-29, 44-34, 44-41, 45-9, 45-12, 45-20, 45-26, 45-30, 45-37, 45-39, 46-11, 46-14, 46-21, 46-37, 46-42, 47-8, 47-13, 47-19, 47-24, 47-31, 47-34, 48-8, 48-11, 48-21, 48-24, 48-31, 48-46, 48-47, 49-7, 49-18, 49-19, 49-22, 49-33, 49-34, 49-40, 50-11, 50-27, 50-31, 50-42, 50-46, 50-48, 51-17, 51-23, 51-32, 51-36, 51-39, 51-41, 52-14, 52-28, 52-30, 52-37, 52-42, 52-46, 53-18, 53-23, 53-33, 53-34, 53-44, 53-49, 54-15, 54-25, 54-28, 54-35, 54-38, 54-43, 55-13, 55-24, 55-29, 55-34, 55-44, 55-47, 56-12, 56-26, 56-29, 56-39, 56-41, 56-45, 57-16, 57-22, 57-27, 57-35, 57-40, 57-43, 58-10, 58-32, 58-33, 58-38, 58-40, 58-43, 59-21, 59-24, 59-26, 59-37, 59-46, 59-48, 60-20, 60-25, 60-30, 60-36, 60-39, 60-45, 61-19, 61-22, 61-31, 61-34, 61-47, 61-49, 62-25, 62-28, 62-30, 62-52, 62-54, 62-60, 63-22, 63-27, 63-31, 63-50, 63-57, 63-61, 64-24, 64-26, 64-29, 64-55, 64-56, 64-59, 65-23, 65-32, 65-33, 65-51, 65-53, 65-58, 66-29, 66-41, 66-44, 66-55, 66-56, 66-64, 67-31, 67-47, 67-48, 67-50, 67-61, 67-63, 68-33, 68-40, 68-49, 68-53, 68-58, 68-65, 69-24, 69-47, 69-48, 69-55, 69-59, 69-64, 69-67, 70-26, 70-37, 70-45, 70-56, 70-59, 70-64, 71-27, 71-35, 71-42, 71-50, 71-57, 71-63, 72-23, 72-41, 72-44, 72-51, 72-53, 72-65, 72-66, 73-28, 73-35, 73-42, 73-52, 73-54, 73-62, 73-71, 74-22, 74-40, 74-49, 74-57, 74-61, 74-63, 74-68, 75-25, 75-36, 75-38, 75-54, 75-60, 75-62, 76-32, 76-36, 76-38, 76-51, 76-58, 76-65, 76-75, 77-30, 77-37, 77-45, 77-52, 77-60, 77-62, 77-70, 78-38, 78-43, 78-54, 78-58, 78-75, 78-76, 79-42, 79-46, 79-50, 79-52, 79-71, 79-73, 80-39, 80-45, 80-56, 80-60, 80-70, 80-77, 81-34, 81-47, 81-55, 81-61, 81-67, 81-69, 82-40, 82-43, 82-57, 82-58, 82-68, 82-74, 82-78, 83-39, 83-41, 83-51, 83-56, 83-66, 83-72, 83-80, 84-34, 84-44, 84-53, 84-55, 84-66, 84-72, 84-81, 85-35, 85-43, 85-54, 85-57, 85-71, 85-73, 85-78, 85-82, 86-34, 86-49, 86-53, 86-61, 86-68, 86-74, 86-81, 86-84, 87-36, 87-39, 87-51, 87-60, 87-75, 87-76, 87-80, 87-83, 88-46, 88-48, 88-50, 88-59, 88-67, 88-69, 88-79, 89-37, 89-46, 89-52, 89-59, 89-70, 89-77, 89-79, 89-88, 90-57, 90-63, 90-71, 90-74, 90-82, 90-85, 91-51, 91-65, 91-72, 91-76, 91-83, 91-87, 92-59, 92-64, 92-69, 92-70, 92-88, 92-89, 93-54, 93-62, 93-73, 93-75, 93-78, 93-85, 94-53, 94-65, 94-68, 94-72, 94-84, 94-86, 94-91, 95-55, 95-64, 95-66, 95-69, 95-81, 95-84, 95-92, 96-56, 96-64, 96-66, 96-70, 96-80, 96-83, 96-92, 96-95, 97-50, 97-63, 97-67, 97-71, 97-79, 97-88, 97-90, 98-52, 98-62, 98-73, 98-77, 98-79, 98-89, 98-93, 99-60, 99-62, 99-75, 99-77, 99-80, 99-87, 99-93, 99-98, 100-61, 100-63, 100-67, 100-74, 100-81, 100-86, 100-90, 100-97, 101-58, 101-65, 101-68, 101-76, 101-78, 101-82, 101-91, 101-94, 102-81, 102-84, 102-86, 102-94, 102-95, 102-100, 103-93, 103-98, 103-99, 104-71, 104-73, 104-79, 104-85, 104-90, 104-93, 104-97, 104-98, 104-103, 105-80, 105-83, 105-87, 105-91, 105-96, 105-99, 105-103, 106-92, 106-95, 106-96, 106-102, 106-103, 106-105, 107-66, 107-72, 107-83, 107-84, 107-91, 107-94, 107-95, 107-96, 107-102, 107-105, 107-106, 108-67, 108-69, 108-81, 108-88, 108-92, 108-95, 108-97, 108-100, 108-102, 108-106, 109-70, 109-77, 109-80, 109-89, 109-92, 109-96, 109-98, 109-99, 109-103, 109-105, 109-106, 110-78, 110-82, 110-85, 110-90, 110-93, 110-101, 110-103, 110-104, 111-91, 111-94, 111-101, 111-102, 111-103, 111-105, 111-106, 111-107, 111-110, 112-75, 112-76, 112-78, 112-87, 112-91, 112-93, 112-99, 112-101, 112-103, 112-105, 112-110, 112-111, 113-79, 113-88, 113-89, 113-92, 113-97, 113-98, 113-103, 113-104, 113-106, 113-108, 113-109, 114-90, 114-97, 114-100, 114-102, 114-103, 114-104, 114-106, 114-108, 114-110, 114-111, 114-113, 115-68, 115-74, 115-82, 115-86, 115-90, 115-94, 115-100, 115-101, 115-102, 115-110, 115-111, 115-114.
\end{spacing}
}

\noindent Eryk Kopczy\'nski, Institute of Informatics, Warsaw University,\\
Banacha 2, 02-097 Warsaw, Poland \\
\texttt{erykk@mimuw.edu.pl}
\\
\\
Igor Pak, Department of Mathematics, UCLA, \\ Los Angeles, CA 90095, USA \\
\texttt{pak(@math).ucla.edu}
\\
\\
Piotr Przytycki, Institute of Mathematics, Polish Academy of Sciences, \\
\'Sniadeckich 8, 00-956 Warsaw, Poland \\
\texttt{pprzytyc@mimuw.edu.pl}

\end{document}